\documentstyle[12pt,amssymb,amscd,psfig]{amsart}
\newtheorem{theorem}{Theorem}
\newtheorem{lemma}[theorem]{Lemma}
\newtheorem{proposition}[theorem]{Proposition}
\newcommand{\co}{\colon\thinspace}

\title[$4$-manifolds with free fundamental group]{Surgery on closed  
$\mathbf{4}$-manifolds with\\ free fundamental group}
\author[V. Krushkal and R. Lee]{Vyacheslav S. Krushkal and Ronnie Lee}
\thanks{VK partially supported by NSF grant DMS 00-72722 and by the Institute for Advanced Study 
(NSF grant DMS 97-29992)} 
\address{Department of Mathematics, Yale University, New Haven, CT 06520-8283, USA}
\email{krushkal\char 64 math.yale.edu, rlee\char 64 math.yale.edu}

\begin{document}

\maketitle


The $4$-dimensional topological surgery conjecture has been established
for a class of groups, including the groups of subexponential growth (see
\cite{FT}, \cite{KQ} for recent developments), however the general case remains
open. The full surgery conjecture is known to be equivalent to 
the question for a class of canonical problems with free fundamental
group \cite[Chapter 12]{FQ}. The proof of the conjecture for ``good'' groups relies 
on the disk embedding theorem (see \cite{FQ}), which is 
not presently known to hold for arbitrary groups. However, in certain
cases it may be shown that surgery works even when the disk embedding
theorem is not available for a given fundamental group (such results still
use the disk-embedding theorem in the simply-connected
setting, proved in \cite{F}.)
For example, this may be done when the surgery kernel is represented by ${\pi}_1$-null
spheres \cite{F1}, or more generally by a ${\pi}_1$-null submanifold satisfying
a certain condition on Dwyer's filtration on second homology \cite{FT2}.
Here we state another instance when the surgery conjecture holds  for free groups.
The following results are stated in the topological category.

\vspace{.1cm}

\begin{theorem} \label{surgery} \sl
Let $X$ be a $4$-dimensional Poincar\'{e} complex with free fundamental group, and 
assume the intersection form on $X$ is extended from the integers.
Let $f\co M\longrightarrow X$ be a degree $1$ normal map, where $M$ is a closed
$4$-manifold. Then the vanishing of the Wall obstruction implies that $f$ is normally 
bordant to a homotopy equivalence $f'\co M'\longrightarrow X$.
\end{theorem}

\vspace{.1cm}

In the canonical surgery problems, $X$ has free fundamental group and trivial ${\pi}_2$,
however what makes them harder to analyze is the interplay between the homotopy type
of $X$, and the topology of the boundary. Our result sidesteps this by
considering closed manifolds. We also prove a related splitting result:

\vspace{.1cm}

\begin{theorem} \label{splitting} \sl
Let $M$ be a closed orientable $4$-manifold with free fundamental group,
and suppose the intersection form on $M$ is extended from the
integers. Then $M$ is s-cobordant to a connected sum of
$\sharp^n S^1\times S^3$ with a simply-connected $4$-manifold.
\end{theorem}

\vspace{.1cm}

Note that if the surgery conjecture fails for free groups, then for both theorems 
above there is, in general, no extension to $4$-manifolds with boundary. 
The assumption on the intersection pairing in theorem \ref{splitting} is necessary,
since there are forms not extended from the integers, for example for ${\pi}_1\cong
{\mathbb{Z}}$ \cite{HT}.
It follows from the classification of $4$-manifolds with infinite cyclic fundamental
group that under the assumptions of theorem \ref{splitting}, if ${\pi}_1(M)\cong{\mathbb{Z}}$
then $M$ is {\em homeomorphic} to a connected sum of $S^1\times S^3$ with a simply-connected
$4$-manifold \cite{FQ}. The $s$-cobordism conjecture for free non-abelian groups
remains open.

\vspace{.2cm}

A brief outline of the proof of theorem \ref{surgery} is as follows.
It may be assumed that ${\pi}_1(f)$ is an isomorphism, and $ker({\pi}_2(f))$ is a 
direct sum of standard planes.
Using the assumption on the intersection form, we construct a complex 
$K=\vee^n S^1\vee K_0$, where $K_0$ is simply-connected, and a map 
$X\longrightarrow K$, inducing isomorphisms on ${\pi}_1$ and ${\pi}_2$.
The inverse preimage of a collection of points $\{ p_i\}$, one in each circle summand
of $K$, under the composition $h\co M\longrightarrow X\longrightarrow K$, is
arranged to be a disjoint union of $3$-spheres in $M$, at the expense of further 
stabilizing $M$. Now consider $2$-spheres in $M$, representing
a hyperbolic basis of $ker ({\pi}_2(f))$; they are surgered along disks lying in
the $3$-spheres, and we show that the resulting elements of ${\pi}_2$ also form a 
hyperbolic basis of the surgery kernel. However, the new $2$-spheres lie in a 
simply-connected $4$-manifold, thus the disk embedding theorem yields
embedded transverse pairs of spheres, concluding the argument.
In the proof of theorem \ref{splitting} one also has to 
keep track of the Lagrangians in the surgery kernel, so that the 
constructed cobordism is a ${\mathbb{Z}}{\pi}_1$-homology product.

Theorem \ref{splitting} also follows from the results of \cite{CH}, \cite{CHR}. However,
there the authors additionally assume that $M$ is smooth, while the conclusion
is still topological. Under the assumption ${\pi}_2(X)=0$, the splitting theorem 
\ref{splitting} is also stated in \cite{H}. 
Our result is entirely in the topological category, and the line of argument is different
from the above papers. In particular, instead of using
$5$-dimensional surgery theory, our proof gives a more explicit geometric construction
of the $s$-cobordism. 

\vspace{.2cm}

\noindent
{\bf Remark.} The idea of our proof extends in certain cases to $4$-manifolds 
with boundary. Recall from the outline above that the proof considers 
point inverses $h^{-1}(p_i)$, changes them into $3$-spheres, up to a 
cobordism of $M$, and reduces the problem to the simply-connected setting. 
For $M$ with non-empty boundary, if there is a map $h\co M\longrightarrow K$
with $h^{-1}(p_i)\cap \partial M=S^2$, for each $i$, then the same proof 
yields $3$-disks, and the argument goes through. Compare this with the general case: 
for example if $h^{-1}(p_i)\cap \partial M$ are tori, they cannot necessarily
be arranged to bound disjoint solid tori, even up to an $s$-cobordism of $M$, see \cite{K}. 
This illustrates the difference between the closed case, considered here, and the canonical 
surgery problems.

\vspace{.2cm}

Before proceeding with the proof of the theorems, we state a preliminary result.
Here we introduce a $2$-complex $K$ which will serve as a reference for the
homotopy data while the cobordisms are being constructed.

\vspace{.2cm}

\begin{lemma}  \label{homotopy} \sl
Let $X$ be a $4$-dimensional Poincar\'{e} complex with free fundamental group,
and suppose the intersection form on $X$ is extended from the
integers. Then there is a simply-connected complex $K_0$ and a map 
$g\co X\longrightarrow K=\vee^n S^1\vee K_0$, inducing isomorphisms 
on ${\pi}_1$ and ${\pi}_2$.
\end{lemma}

\vspace{.2cm}

\noindent
{\em Proof.} 
By \cite[\S 0.2]{W} we may replace $X$ by a homotopy equivalent $4$-complex 
with a single top cell. Following \cite{MWh}, consider the
$3$-type $({\pi}_1, {\pi}_2, k)$ of its $3$-skeleton $X^{(3)}$. Here ${\pi}_1=
{\pi}_1(X)$, ${\pi}_2={\pi}_2(X)$ is a module over ${\pi}_1$, and the invariant 
$k$ vanishes, since $k$ is an element of $H^3({\pi}_1; {\pi}_2)=0$. The assumption 
on the intersection form implies that ${\pi}_2(X)\cong A\otimes_{\mathbb{Z}}
{\mathbb{Z}}{\pi}_1$ where $A$ is an abelian group, and the intersection form
on ${\pi}_2(X)$ is induced from a form
 on $A$.
Let $K_0=K(A,2)$, and set 

\[ K=K({\pi}_1,1)\vee K(A,2) =\vee^n S^1\vee K_0,\]

\noindent
then $K^{(3)}$ has the same algebraic $3$-type as 
$X^{(3)}$. It is proved in \cite{MWh} that any homomorphism of $3$-types is 
induced by a map between the $3$-complexes. Hence there is a map 
$g\co X^{(3)}\longrightarrow K$, inducing isomorphisms on ${\pi}_1$ and ${\pi}_2$.

We claim that the obstruction $o$ to extending $g$ over the top cell of $X$, $o\in
H^4(X;{\pi}_3(K))$, vanishes.
Since ${\pi}_3(K_0)=0$, by Hilton-Milnor's theorem ${\pi}_3(K)$ is generated by 
the Whitehead products $[{\alpha}\cdot a, {\beta}\cdot b]$, where $a,b\in{\pi}_2(K_0)$, 
and ${\alpha}\neq{\beta}\in{\pi}_1$. Suppose $o$ does not vanish, then the value 
of $o$ on the fundamental  cycle in $X$ is non-trivial in ${\pi}_3(K)$. Since $g$ 
induces isomorphisms on ${\pi}_1$ and ${\pi}_2$, the attaching map of the $4$-cell
has a non-trivial component $[{\alpha}\cdot a, {\beta}\cdot b]\in{\pi}_3(X^{(3)})$, where
${\alpha}\neq {\beta}\in {\pi}_1$, and $a,b,\in A$.
However, the intersection number of two classes in ${\pi}_2(X)$ is determined
by the value of the attaching map of the $4$-cell on their Whitehead product
in $X^{(3)}$. In particular, in the situation above the intersection
$({\alpha} a)\cdot ({\beta} b)$ is in ${\mathbb{Z}}$,
contradicting the assumption on the intersection pairing on $X$, and the
assumption ${\alpha}\neq{\beta}$.
\qed

\vspace{.5cm}

\noindent
{\em Proof of theorem} \ref{surgery}. 
Following the proof on the higher-dimensional surgery theorem \cite{W} (see also
\cite[Chapter 11]{FQ}), we may assume that $f$ induces an isomorphism on ${\pi}_1$, 
and the kernel of ${\pi}_2(f)$ is a direct sum of standard planes. Since
the intersection pairing on $X$ is induced from the integers, the same is true for
the intersection form on $M$.
Consider a map $g\co X\longrightarrow K$ given by lemma \ref{homotopy}, and arrange
the composition $h=gf\co M\longrightarrow K$ to be transverse to 
a collection of points $p_1,\ldots, p_n$, different from the basepoint,
one in each circle summand of $K$. (See \cite{Q} or \cite[\S 9.6]{FQ} for the statement
of transversality in the topological category.) Denote the $3$-manifold
$h^{-1}(p_i)$ by $P_i$, and set $P=\amalg^n P_i$. Changing the map $h$ by a homotopy
if necessary, we may assume that $P_i$ is connected, for each $i$.

For $i=1,\ldots, n$, consider a framed link $L_i\subset P_i$ such that the surgery on 
$P_i$ along $L_i$ gives the $3$-sphere (cf \cite{Ka}), and let $L=\amalg L_i$. 
If the components of $L$ bounded disjoint embedded disks with interiors
in $M\smallsetminus P$ and with appropriate framings, then $P$ could be ambiently 
surgered to get a collection of disjoint $3$-spheres, geometrically dual to the 
generators of ${\pi}_1(M)$. Since this cannot be expected in general, we perform 
surgery along the link $L$ on the $4$-manifold $M$, and denote the result by $N$.  
Here for the surgeries on $M$ we use the framing of $L$, determined by the
framing of $L$ in $P$. The components of $L$ bound disjoint embedded disks
with the required framings in $N$, thus the map $h$ is bordant
to $h'\co N\longrightarrow K$, with $(h')^{-1}(p_i)=S^3$.

Since $L$ is null-homotopic in $M$, its components bound disjoint embedded disks 
in $M$. This implies that $N$ is homeomorphic to $M$, connected summed
with several copies of $S^2\times S^2$, and also with copies of the twisted bundle
$S^2\widetilde{\times} S^2$. The framed link $L$ may be chosen so that 
$P\times I\cup_L 2$-handles is spin \cite{Ka}, so if ${\omega}_2(M)=0$ then
${\omega}_2(N)$ is also trivial. 
If ${\omega}_2(M)\neq 0$, note that $H_2(M;{\mathbb{Z}})$ consists of spherical
classes, thus $M\sharp S^2\widetilde{\times} S^2\cong M\sharp S^2\times S^2$.
In either case, we may assume $N\cong M\sharp^k(S^2\times S^2)$, and we have
a normal bordism from the map $f\co M\longrightarrow X$ to a map $N\longrightarrow X$.

At this point we caution that without a restriction on the intersection pairing on $M$, one 
could still consider a map $M\longrightarrow \vee^n S^1$, classifying ${\pi}_1$, and construct 
$N$ as above. (A similar construction is used in \cite{H} in a proof of the stable $4$-dimensional
Kneser's conjecture.) Now the intersection form on $N$ (the form on $M$, stabilized by adding 
several hyperbolic pairs), is extended from the integers, since $N$ contains $3$-spheres,
geometrically dual to the generators of ${\pi}_1$. Certainly, in this general case one
cannot hope to de-stabilize $N$ while preserving the $3$-spheres.

Returning to the proof, denote $N_{0}={\rm closure}(N\smallsetminus (\amalg S^3\times I))$.
Then ${\pi}_2(N)\cong {\pi}_2(N_{0})\otimes_{\mathbb{Z}} 
{\mathbb{Z}}{\pi}_1$, and moreover 

\[ker[{\pi}_2(N)\longrightarrow {\pi}_2(K)]
\cong  ker[{\pi}_2(N_{0})\longrightarrow
{\pi}_2(K_0)]\otimes_{\mathbb{Z}}{\mathbb{Z}}{\pi}_1. \]

\noindent
Consider a standard hyperbolic basis for $ker[{\pi}_2(N)\longrightarrow {\pi}_2(K)]$,
say $\{a_i, b_i\}$. Arrange these $2$-spheres to be 
transverse to $\amalg S^3$, and let the circles of intersection bound maps of disks 
in the $3$-spheres. Consider two copies of each disk, lying in $S^3\times \{ -{\epsilon}\}$
and $S^3\times \{ {\epsilon} \} $ respectively, and use them to surger the spheres $a_i, b_i$.
In other words, we cut out an annulus out of each $2$-sphere, and glue in the disks described above.
The resulting $2$-spheres lie in the complement of $\amalg S^3$, and we connect them
to the basepoint by arcs in $N_0$. 
The constructed classes ${\alpha}_i, {\beta}_i
\in {\pi}_2(N_{0})$ are homologous, but not necessarily homotopic to $a_i, b_i$.
For example, suppose $b_i$ intersects $S^3$ in a circle. Cutting $b_i$ as above, we get
two spheres $b'_i$ and $b''_i$ with $b_i=b'_i+b''_i\in {\pi}_2(N)$, while ${\beta}_i=b'_i+
g \, b''_i$. Here $g$ is the generator of ${\pi}_1(N)$ dual to the given $S^3$. 

Since $\{ {\alpha}_i, {\beta}_i\}$ are homologous to the original hyperbolic  basis, these classes 
freely generate 

\[ ker[H_2(N_0)\longrightarrow H_2(K_0)]\cong ker[{\pi}_2(N_0)\longrightarrow {\pi}_2(K_0)], \]

\noindent 
hence they also freely generate 

\[ ker[{\pi}_2(N)\longrightarrow {\pi}_2(K)] =  ker[{\pi}_2(N)\longrightarrow 
{\pi}_2(X)] \]

\noindent
as a module over ${\pi}_1$. Moreover, $\{ {\alpha}_i, {\beta}_j\}$ is a collection
of (algebraically) transverse pairs of spheres in $N_0$, and the disk embedding theorem
in the simply-connected setting \cite{F}, \cite[\S 5.1]{FQ} gives
a collection of embedded transverse pairs, homotopic to $\{ {\alpha}_i,
{\beta}_i\}$. Surgering them out yields a homotopy equivalence $f'\co M'\longrightarrow X$.
\qed

\vspace{.5cm}

\noindent
{\em Proof of theorem} \ref{splitting}. Since $M$ is homotopy equivalent to a
Poincar\'{e} complex (cf \cite[Chapter 3]{KS}), there is a map $f\co M\longrightarrow K$, 
satisfying the conclusions of lemma \ref{homotopy}. As in the proof of theorem
\ref{surgery}, arrange it to be transverse to a collection of points $p_1,\ldots, p_n$, 
one in each circle summand of $K$, and denote $P_i=f^{-1}(p_i)$, $P=\amalg^n P_i$. 
We may assume that $P_i$ is connected, for each $i$. 

For $i=1,\ldots, n$, consider a framed link $L_i\subset P_i$ such that the surgery on 
$P_i$ along $L_i$ gives the $3$-sphere, and let $L=\amalg L_i$. We denote the surgery 
on $M^4$ along the link $L$ with the corresponding framings by $N$. 
Define 

\[ W_1=M\times [0,1]\cup 2 {\rm -handles}; \hspace{.1cm} \partial W_1=M\amalg N. \]

\noindent
As in the proof of theorem \ref{surgery}, the link $L$ may be chosen so that 
$N\cong M\sharp^k S^2\times S^2$.
The intersection form over ${\mathbb{Z}}{\pi}_1$ on $N$ is the form on $M$, plus $k$
standard planes. We fix notation, $\{ a_i, b_i \}$,  for a hyperbolic basis of
$ker[{\pi}_2(N)\longrightarrow{\pi}_2(K)]$, where $\{ b_i\}$ correspond to the belt spheres
of the $2$-handles of $W_1$. Let $A,B$ be the 
${\mathbb{Z}}{\pi}_1$-submodules of ${\pi}_2(N)$, (freely) generated by the 
$\{ a_i\}$ and $\{ b_i\}$ respectively, then 

\[ ker[{\pi}_2(N)\longrightarrow{\pi}_2(K)]\cong A\oplus B. \]

\noindent
Note that the homomorphism ${\phi}\co A\longrightarrow B$, induced by the intersection
pairing: ${\phi}(a)={\Sigma}_j(a\cdot b_j) b_j$, is an isomorphism: ${\phi}(a_i)=b_i$ for each $i$.

The components of $L$ bound in $N$ disjoint embedded framed disks, provided by 
the $4$-dimensional surgeries. Use these disks to ambiently surger each 
$P_i=f^{-1}(p_i)$ in $N$ into the $3$-sphere. 
Denote $N_{0}={\rm closure}(N\smallsetminus(\amalg S^3\times I))$.
As in the proof of theorem \ref{surgery}, consider $2$-spheres ${\alpha}_i$, ${\beta}_i
\subset N_{0}$, homologous (but not necessarily homotopic) to $a_i$, $b_i$.
These are obtained by arranging the intersections $a_i, b_i\cap\amalg S^3$ to be transverse;
each circle of intersection bounds a map of a disk in $S^3$, and we surger $a_i$ and $b_i$ 
along these disks. Finally, connect the resulting spheres to the basepoint by arcs
in $N_{0}$ to get ${\alpha}_i$, ${\beta}_i$.

\vspace{.3cm}

\begin{proposition} \label{basis} \sl
The classes $\{ {\alpha}_i \}$ freely generate the ${\mathbb{Z}}{\pi}_1$-module $A$.
\end{proposition}

\vspace{.3cm}

\noindent
{\em Proof.} 
The map $f\co M\longrightarrow K$ extends to a map $g\co N\longrightarrow K\vee^k 
(S^2\times S^2)$, inducing isomorphisms on ${\pi}_1$ and ${\pi}_2$, and so that 
${\pi}_2(g)$ maps $a_i,b_i$ to the generators $\bar a_i, \bar b_i$ of the corresponding
${\pi}_2(S^2\times S^2)$. Here $k$ is the number of components of the link $L$. Set 

\[ K'_0=K_0\vee^k (S^2\times S^2)\cup 3{\rm -cells}, \hspace{.2cm} K'=\vee^n S^1\vee K'_0 \] 

\noindent
where the 3-cells are attached to $\{ \bar a_i \}$, and note that $A=ker[{\pi}_2(N)
\longrightarrow {\pi}_2(K')]$. Since $\{ {\alpha}_i\}$ are homologous to $\{ a_i\}$, 
the collection $\{ {\alpha}_i\}$ freely generates 

\[ ker[H_2(N_{0})\longrightarrow H_2(K'_0)]
\cong ker[{\pi}_2(N_{0})\longrightarrow {\pi}_2(K'_0)]; \] 

\noindent 
thus they are also free generators over ${\mathbb{Z}}{\pi}_1$ of 

\[ ker[{\pi}_2(N_{0})\longrightarrow {\pi}_2(K'_0)]\otimes_{\mathbb{Z}} {\mathbb{Z}}{\pi}_1
\cong ker[{\pi}_2(N)\longrightarrow {\pi}_2(K')]\cong A. \qed  \]

The spheres $\{ {\alpha}_i$, ${\beta}_j\}$ form a collection of algebraically transverse 
pairs in $N_0$, and the disk embedding theorem in the simply-connected setting 
\cite{F}, \cite[\S 5.1]{FQ} implies that they are homotopic to embedded transverse
pairs ${\alpha}'_i$, ${\beta}'_i$. Set 

\[ W_2=N\times [0,1]\cup 3{\rm -handles} \] 

\noindent where the $3$-handles are attached to the spheres ${\alpha}'_i\subset
N\times\{ 1\}$, and let $M'$ be the corresponding surgery on $N$. Consider the 
cobordism $W=W_1\cup_N W_2$ between $M$ and $M'$. Since there are only $2$- 
and $3$-handles, the chain complex for the relative ${\mathbb{Z}}{\pi}_1$ homology groups
is $0\longrightarrow C_3\longrightarrow C_2\longrightarrow 0$, where $C_i$ is the
${\mathbb{Z}}{\pi}_1$-module, freely generated by the $i$-handles, cf \cite{M}. 
The boundary homomorphism is given by the intersection numbers of the 
attaching spheres of the $3$-handles with the belt spheres of the $2$-handles. 
Using proposition \ref{basis}, observe that the homomorphism $C_3\longrightarrow C_2$ is 
identified with the isomorhism ${\phi}\co A\longrightarrow B$ considered earlier in the 
proof. Thus $W$ is an $h$-cobordism, and since ${\pi}_1(M)$ is free, its Whitehead group
is trivial, and so $W$ is an $s$-cobordism as asserted in the theorem.
\qed

\end{document}